\documentclass[11pt]{elsarticle}
\usepackage[utf8]{inputenc}
\usepackage{amsmath,amssymb,amsfonts,amsthm}
\usepackage{bm}
\usepackage[margin=3cm]{geometry}
\usepackage{quoting,xparse}
\usepackage{pgfplots}
\pgfplotsset{compat=newest}
\pgfplotsset{plot coordinates/math parser=false}

\usepackage{abstract}
\usepackage{hyperref}
\usepackage[capitalize,nameinlink]{cleveref}

\usepackage{algorithm}
\usepackage{algpseudocode}

\usepackage{dsfont}
\usepackage{todonotes}
\usepackage{verbatim}
\usepackage{multicol}
\usepackage{enumitem}
\usepackage{graphicx} 
\usepackage[english]{babel}
\usepackage{mathtools}
\newtheorem{definition}{Definition}
\graphicspath{ {images/} }
\allowdisplaybreaks
\newtheorem{theo}{Theorem}[section]

\newtheorem{cor}[theo]{Corollary}

\newtheorem{rem}[theo]{Remark}
\usepackage{mathtools}

\newenvironment{talign*}
 {\csname align*\endcsname}
 {\endalign}

 \newenvironment{tgather*}
 {\csname gather*\endcsname}
 {\endgather}

\title{Hybrid multi-population traffic flow model: Optimal control for a mean-field limit\\}

\author[1]{Maria Teresa Chiri}
\author[2]{Christopher Denaro}
\author[3]{Xiaoqian Gong}
\author[4]{Benedetto Piccoli}

\address[1]{Department of Mathematics and Statistics, Queen's University,
Kingston, ON 
Canada, {\tt\small maria.chiri@queensu.ca}}
\address[2]{Department of Mathematical Sciences and Center for Computational and Integrative Biology, Rutgers University, Camden, NJ, USA, {\tt\small cad373@scarletmail.rutgers.edu}}
\address[3]{Department of Mathematics, Amherst College
Amherst, MA, USA, {\tt\small xgong@amherst.edu}}
\address[4]{Department of Mathematical Sciences and Center for Computational and Integrative Biology, Rutgers University, Camden, NJ, USA, {\tt\small piccoli@camden.rutgers.edu}}

\begin{document}

\maketitle

\begin{abstract}
Modeling heterogeneous and multi-lane traffic flow is essential for understanding and controlling complex transportation systems. In this work, we consider three vehicle populations: two classes of human-driven vehicles (cars and trucks) and autonomous vehicles, the latter characterized by controlled acceleration.
Compared to single-population models, multi-population modeling poses greater challenges, primarily due to the increased number of parameters required to describe lane-changing behavior and the added complexity in passing to the mean-field limit. 
We model multi-lane traffic as a hybrid dynamical system, combining continuous dynamics within each lane and discrete events corresponding to lane-changing maneuvers. We then formulate and analyze the optimal control problem associated with such hybrid systems from both microscopic and mesoscopic perspectives. Using techniques from $\Gamma$-convergence, we prove the existence of solutions to the optimal control problem in the mean-field limit of a finite-dimensional hybrid system. Finally, we present numerical simulations illustrating the impact of trucks on overall traffic efficiency.
\end{abstract}

\section{Introduction}
\label{sec:introduction}
A key challenge in traffic flow modeling is the representation of multi-lane traffic. The complexity arises from its hybrid nature: each lane exhibits continuous dynamics, while lane changes introduce discrete events. Lane changing is a frequent yet high-risk maneuver involving multiple vehicles  \cite{HertyVisconti2018}. Moreover, an accurate description must account for the diversity of vehicle types, such as cars and trucks, which may exhibit different lane-changing behaviors.
Currently, multi-lane traffic is modeled either using two-dimensional approaches~\cite{HertyMoutariVisconti2018, SukhinovaTrapeznikovaChetverushkinChurbanova2009}, where lane-changing rules are not explicitly defined, or by treating lanes as discrete entities~\cite{HoldenRisebro2019,SongKarni2019}. The heterogeneous composition of traffic flow is becoming increasingly important, particularly with the growing interest in automated vehicles and their impact on vehicular traffic~\cite{HoogendoormReview2014}. Both experiments~\cite{DelleMonache2019, Piccoli-DissipationStopAndGo2018} and mathematical models~\cite{TosinZanella2021,Piccoli2020} have demonstrated that even a small number of controlled vehicles can help stabilize traffic flow and mitigate unstable phenomena.

In \cite{gpv_meanfield}, the authors introduced mean-field equations coupled with ODEs to model the dynamics of mixed traffic, where human-driven vehicles are represented by a mean-field equation and a small number of autonomous vehicles by ODEs. The system is formally derived from a microscopic model that combines the Bando model \cite{bando1994structure} with the Follow-the-Leader (FtL) model \cite{reuschel1950vehicle1}, along with lane-changing conditions inspired by \cite{KTH07}. The case of two populations—human-driven vehicles (cars and trucks) and autonomous vehicles—was studied in \cite{CGPmultipopulation}.
Since trucks are larger and slower to accelerate to their desired speeds than cars, they exert a greater influence on traffic flow. Consequently, an accurate traffic model should not overlook the impact of this population.

In this work, we generalize previous results by incorporating cars, trucks, and autonomous vehicles while addressing a general optimal control problem. This framework applies to the optimization of travel time, fuel consumption, and emissions. Due to the hybrid nature of the dynamics, the problem falls into the category of hybrid optimal control problems, as studied in \cite{MR1786332, 827981}.

Specifically, the control variable corresponds to the regulated acceleration of autonomous vehicles, which, in turn, influences the overall dynamics of the other two populations. The cost function is formulated in a general form, incorporating terms that depend on the position, velocity, control, and lane-switching times of all vehicles. Consequently, our approach can be used for tracking a desirable traffic profile or minimizing travel time, fuel consumption, and CO2 emissions \cite{Piccoli-DissipationStopAndGo2018, YAO2021120766}.

A vast body of literature exists on the optimal control of hybrid systems, though most results focus on finite-dimensional dynamics \cite{654885, garavello2005hybrid,piccoli1998hybrid}. Here, we first formulate the optimal control problem for the finite-dimensional (microscopic) setting with three interacting populations. We then pass to the mean-field limit and prove $\Gamma$-convergence \cite{DalMaso} to the optimal control problem for the coupled system: a mean-field partial differential equation governing car and truck traffic and a system of controlled ODEs for autonomous vehicles. This approach ensures both the existence of a solution to the limiting problem and the convergence of optimal control policies.

The paper is organized as follows. In Section \ref{Prel}, we introduce a microscopic model for heterogeneous traffic, incorporating cars, trucks, and autonomous vehicles on $L$ lanes. We then specify the lane-changing dynamics in terms of acceleration and a probability measure that captures the uncertainty in lane-changing decisions.

To prevent Zeno phenomena \cite{1656623} and other singularities, we impose a cool-down time assumption: consecutive lane-changing maneuvers are not allowed without a minimum time separation. While this assumption excludes some extreme behaviors, it remains consistent with data collected on highways \cite{HertyVisconti2018}.

In Section \ref{opt_FD}, we formally define the controlled finite-dimensional hybrid system and its associated optimal control problem. The mean-field limit is presented in Section \ref{inf_dim}, where we also prove the $\Gamma$-convergence of the optimal control problem.

Section \ref{gener} provides a preliminary extension of the hybrid system to $M+1$ vehicle populations, illustrating the flexibility of our framework in capturing more refined multi-population traffic models. In Section \ref{sec:simulations}, we present numerical simulations based on the finite-dimensional hybrid system, highlighting the impact of the truck penetration rate on overall traffic performance.

Finally, the last section summarizes our contributions, emphasizing the significance of modeling mixed traffic with cars and trucks and outlining future research directions.

\section{Preliminary}\label{Prel}
In this section, we begin by formulating the convolutional version of a second-order car-following model—the Bando-Follow-the-Leader model—tailored for heterogeneous traffic including cars, trucks, and autonomous vehicles. We then describe the lane-changing rules for multi-lane and multi-class traffic, addressing both finite-dimensional and infinite-dimensional frameworks. The model includes parameters that account for the distinctions between vehicle types and the interaction dynamics among them. Finally, we highlight a key assumption of our lane-changing framework: the \emph{cool-down time}, which prevents vehicles from performing consecutive lane changes in rapid succession.

\subsection{Car-following models for heterogeneous traffic }\label{SecA}

For clarity, we introduce the notations used throughout the paper. Let \( T > 0 \) be fixed, and let \( \mathcal{I} \) denote the set of vehicle indices on an open stretch of road of \(L\) lanes. Let \(l_i \in \mathbb{R}_{>0}\) be the length of \(i \in \mathcal{I}\). We define the set of lane labels as  \( \mathcal{K} = \{1, \dots, L\} \). We denote by \( \mathcal{P} \) the space of probability measures and by \( \mathcal{M}_+ \) the space of positive Borel measures. We then use the following notations:
\begin{itemize}
    \item \( \mathcal{I}_P, \mathcal{I}_Q, \mathcal{I}_S \): index sets for cars, autonomous vehicles, and trucks on the road, respectively.
    \item \( \mathcal{I}_P^k, \mathcal{I}_Q^k, \mathcal{I}_S^k \): index sets for cars, autonomous vehicles, and trucks on lane \( k \in \mathcal{K} \), respectively.
    \item \( (x_i(t), v_i(t)) \): position and velocity of vehicle \( i \in \mathcal{I}\) at time \( t \in [0, T] \).
    \item \( P_k(t), Q_k(t), S_k(t) \): number of cars, autonomous vehicles, and trucks, respectively, on lane \( k \in \mathcal{K} \) at time \( t \in [0, T] \).
\end{itemize}

The second-order car-following model, known as the Bando-Follow-the-Leader (Bando-FtL) model, was originally introduced in \cite{Piccoli-DissipationStopAndGo2018}. It assumes that the acceleration of a vehicle \( i \in \mathcal{I} \) depends on its space headway \( h_i \coloneqq x_{L_i} - x_i - l_i \), the relative velocity with respect to its leading vehicle \( \Delta v_i \coloneqq v_{L_i} - v_i \), and the deviation of the velocity of vehicle \(i\), \(v_i\), from its optimal velocity \( V_i \). Here, \( L_i \) denotes the index of the leading vehicle of \( i \in \mathcal{I} \). Specifically, the Bando-FtL model is governed by the following dynamics
\begin{align*}
    \dot{x}_i & = v_i, && i \in \mathcal{I}\\
    \dot{v}_i & = \alpha(V_i - v_i) + \beta \tfrac{\Delta v_i}{h_i^2} && i \in \mathcal{I}\setminus \{i \colon i \text{ is the index of the leading vehicle}\}
    \end{align*}
where \(\alpha, \beta \in \mathbb{R}_{>0}\) are the model parameters.

The following convolution-based formulation of the Bando--Follow-the-Leader (FtL) model for homogeneous traffic consisting of a single class of \( N \) vehicles was introduced in \cite{gpv_meanfield}:
\begin{align*}
	\dot{x}_i & = v_i, && i \in \{1, \dots, N\}, \\
	\dot{v}_i & = \left(H_1 *_1 \mu + H_2 * \mu\right)(x_i, v_i), && i \in \{1, \dots, N\},
\end{align*}
where \(\mu \in \mathcal{P}(\mathbb{R} \times \mathbb{R}_0^+)\) denotes an atomic probability measure supported on the absolutely continuous trajectories \(t \in [0, T] \mapsto (x_i(t), v_i(t))\), \(i \in \{1. \dots, N\}\), that is,
\begin{equation}
\label{eqn_atomic_measure_P}
\mu(t) = \frac{1}{N} \sum_{i=1}^{N} \delta_{(x_i(t), v_i(t))},
\end{equation}
where for any \((x, v) \in \mathbb{R} \times \mathbb{R}_{\geq 0}\), $\delta_{(x,v)} \in \mathcal{P}(\mathbb{R} \times \mathbb{R}_{\geq 0})$ is the Dirac measure concentrated on $(x,v)$.
Here, \(H_1 *_1 \mu\) represents the convolutional form of the Bando interaction term, and \(H_2 * \mu\) corresponds to the convolutional form of the FtL interaction. The functions \(H_1 \colon \mathbb{R} \times \mathbb{R}_0^+ \to \mathbb{R}\) and \(H_2 \colon \mathbb{R} \times \mathbb{R} \to \mathbb{R}\) are the respective convolution kernels. Specifically, {$\ast_1$ represents the convolution with respect to the first variable and} $\ast$ is the convolution with respect to both variables.

In this work, we introduce a convolutional formulation of the Bando-FtL model tailored to heterogeneous traffic comprising cars, trucks, and autonomous vehicles.

\begin{figure}
     \centering
     \includegraphics[scale=0.2]{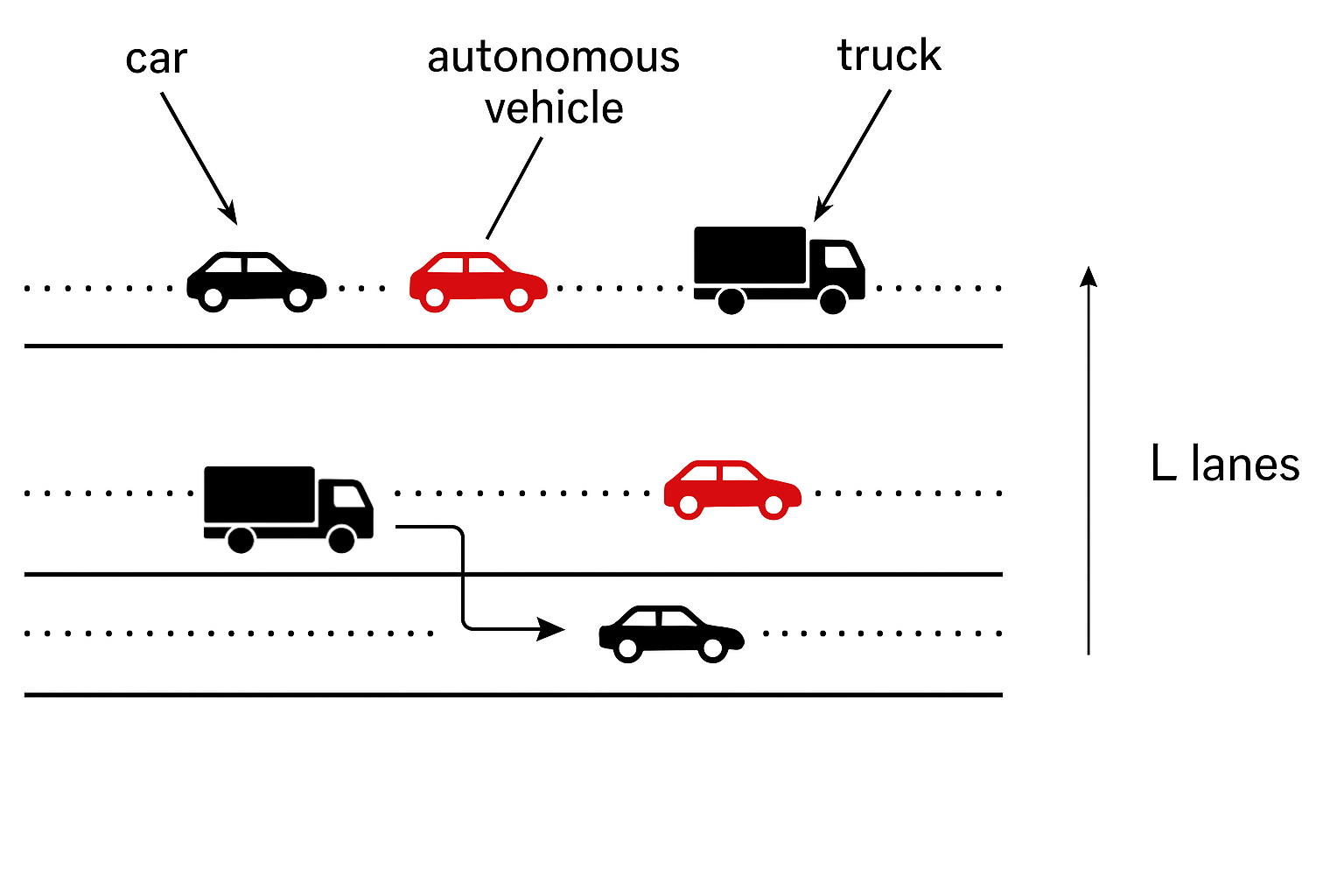}
     \caption{A schematic view of the setting: we consider autonomous vehicles (red), cars (black), and trucks (black) traveling on \(L\) lanes. Vehicles may change lanes according to a probability given by \eqref{eqn: probability_function}, which depends on the Incentive and Safety criteria specified in \eqref{lane-changing_acc}--\eqref{lane-changing_safe}.}
\end{figure}

For each lane \(k \in \mathcal{K}\), we define the atomic measures \(\mu_P^k\), \(\mu_Q^k\), and \(\mu_S^k\), each supported on absolutely continuous trajectories \(t \in [0, T] \mapsto (x_i(t), v_i(t)) \in \mathbb{R} \times \mathbb{R}_{\geq 0}\), for \(i \in \mathcal{I}\), as follows:
\begin{align}
    \mu_P^k(t) &= \frac{1}{P_k(t)} \sum_{i \in \mathcal{I}_P^k(t)} \delta_{(x_i(t), v_i(t))}, \\
    \mu_Q^k(t) &= \frac{1}{Q_k(t)} \sum_{i \in \mathcal{I}_Q^k(t)} \delta_{(x_i(t), v_i(t))}, \label{emAV} \\
    \mu_S^k(t) &= \frac{1}{S_k(t)} \sum_{i \in \mathcal{I}_S^k(t)} \delta_{(x_i(t), v_i(t))}.
\end{align}

These atomic measures encode the positions and velocities of the three vehicle types—cars, trucks, and autonomous vehicles—on each lane. To account for the four possible vehicle-type interactions, Car-Car, Car-Truck, Truck-Car, Truck-Truck, we define the following family of convolution kernels 
\begin{align*}
H_1^{\gamma}: &\mathbb{R}\times\mathbb{R}_{\geq 0 }\rightarrow \mathbb{R} \quad \hbox{ with } \gamma\in\lbrace cc,tc,ct,tt \rbrace\\
& (x,v)\mapsto\alpha_{\gamma} h_{\gamma}(x)(V_{\gamma}(-x)-v)
\end{align*}
where $\alpha_{\gamma}$ are positive parameters denoting the speed of response, $V_{\gamma}$ is the optimal velocity function, and $h_{\gamma}:\mathbb{R}\mapsto\mathbb{R}_{\geq 0 }$ is a smooth function with compact support $[-\varepsilon_{\gamma}, 0]$ measuring the interaction of two vehicles depending on their distance and types, where $\varepsilon_{\gamma} >0$.  A possible choice for $h_{\gamma}$ is the following
\[h_{\gamma}(x) =e^{-\Big(\left(\frac{\varepsilon_{\gamma}}{2}\right)^2 - \left(-x-\frac{\varepsilon_{\gamma}}{2}\right)^2\Big)^{-1}} \cdot \chi_{S}(x),\]
where $A=\{x|-\epsilon_{\gamma} <x <0\}$ and $\chi_{A}$ is the indicator function of set $A$.
Then, for instance, for a car $i \in \mathcal{I}_{P}^k$ on lane $k \in \mathcal{K}$, the Bando-term of the Bando-FtL model can be written as 
\begin{align*}
\left(H_1^{cc}\ast_1 (\mu_P^k+\mu_Q^k)+H_1^{tc}\ast_1\mu_S^k\right)(x_i,v_i), i \in \mathcal{I}_{P}^k.
\end{align*}
Similarly, we introduce the following convolution kernels for the FtL term of the Bando-FtL model, 
\begin{align*}
H_2^{\gamma}: &\mathbb{R}\times\mathbb{R}_{\geq 0 }\rightarrow \mathbb{R} \quad \hbox{ with } \gamma\in\lbrace cc,tc,ct,tt \rbrace\\
& (x,v)\mapsto\beta_{\gamma} h_{\gamma}(x)\frac{-v}{x^2}
\end{align*}
where $\beta_{\gamma}$ is positive. Therefore, for a car $i \in \mathcal{I}_P^k$ on lane $k \in \mathcal{K}$, the FtL-term can be written as 
\begin{align*}
\left(H_2^{cc}\ast (\mu_P^k+\mu_Q^k)+H_2^{tc}\ast\mu_S^k\right)(x_i,v_i), i \in \mathcal{I}_{P}^k.
\end{align*}
Additionally, we assume that the dynamics of autonomous vehicles are distinguished from those of other vehicles primarily through the inclusion of control terms. In particular, the acceleration of each autonomous vehicle \( i \in \mathcal{I}_{Q}^k \) is governed by the Bando--FtL model augmented with a control input \( u_i \).\\
The vehicles' dynamics are governed by the following system: 
\begin{align}
\label{finite_dimensional_system}
    &\dot{x}_i=v_i \quad i \in \mathcal{I}\notag\\
   &\dot{v}_i=\begin{cases} \left(H_1^{cc}\ast_1 (\mu_P^k+\mu_Q^k)+H_1^{tc}\ast_1\mu_S^k\right)(x_i,v_i)\\
 \qquad +\left(H_2^{cc}\ast(\mu_P^k+\mu_Q^k)+H_2^{tc}\ast\mu_S^k\right)(x_i,v_i),
     \qquad i \in \mathcal{I}_P^k,
     \\
     \\
     \left(H_1^{cc}\ast_1 (\mu_P^k+\mu_Q^k)+H_1^{tc}\ast_1\mu_S^k\right)(x_i,v_i)\\
 \qquad +\left(H_2^{cc}\ast(\mu_P^k+\mu_Q^k)+H_2^{tc}\ast\mu_S^k\right)(x_i,v_i) +u_i, 
     i \in \mathcal{I}_Q^k,
     \\
     \\
     \left(H_1^{ct}\ast_1 (\mu_P^k+\mu_Q^k)+H_1^{tt}\ast_1\mu_S^k\right)(x_i,v_i)\\
    \qquad +\left(H_2^{ct}\ast(\mu_P^k+\mu_Q^k)+H_2^{tt}\ast\mu_S^k\right)(x_i,v_i), \qquad \,i \in \mathcal{I}_S^k.
     \end{cases}
\end{align}

\subsection{Lane-changing rules for multi-lane traffic}
In the microscopic (finite-dimensional) setting, lane-changing behavior is governed by safety and incentive conditions. In contrast, the mean field limit considers the number of vehicles tending to infinity; accordingly, we formulate a set of rules that capture the corresponding macroscopic dynamics.\\

\textbf{ The finite-dimensional case}\\

We consider $P$ cars, $Q$ autonomous vehicles and $S$ trucks on the road of \(L\) lanes. Let $a_i^k$ be the acceleration of vehicle $i \in \mathcal{I}$ on lane $k \in \mathcal{K}$. If vehicle $i$ moves to the lane $k' \in \{k-1, k+1\} \cap \mathcal{K}$ , then we denote with  $\bar{a}_i^{k'}$ the expected acceleration on this new lane. Let $i_F^{k'}$ be the index of  the following vehicle of vehicle $i$ on its target lane $k'$ if vehicle $i$ performs lane-changing from lane $k$ to lane $k'$ and $\bar{a}_{i_F^{k'}}^{k'}$ the expected acceleration of vehicle \(i_F^{k'}\) on the new lane.

Let $\Delta^{\gamma}$ be positive constants  with $\gamma\in\lbrace c, t, cc,tc,ct,tt \rbrace$. The  {incentive} and  {safety} conditions for lane-changing are defined as follows: 
\begin{align}
\label{lane-changing_acc}
\text{Incentive: } \bar{a}_i^{k'} \geq
 \begin{cases}
 a_i^k+\Delta^{cc} \hbox{ if } i, i_F^{k'} \in \mathcal{I}_P \cup \mathcal{I}_Q,\\
 a_i^k+\Delta^{tc} \hbox{ if } i \in \mathcal{I}_P \cup \mathcal{I}_Q, i_F^{k'} \in \mathcal{I}_S,\\
 a_i^k+\Delta^{ct} \hbox{ if } i\in \mathcal{I}_S, i_F^{k'} \in \mathcal{I}_P \cup \mathcal{I}_Q,\\
 a_i^k+\Delta^{tt} \hbox{ if } i, i_F^{k'} \in \mathcal{I}_S;\\
\end{cases}
\end{align}
\begin{align}
\label{lane-changing_safe}
\text{Safety: } \bar{a}_i^{k'}\geq \begin{cases}
-\Delta^{c} \text{ and } \bar{a}_{i_F^{k'}}^{k'} \geq -\Delta^{c}\hbox{ if } i, i_F^{k'} \in \mathcal{I}_P \cup \mathcal{I}_Q,\\
-\Delta^{c} \text{ and } \bar{a}_{i_F^{k'}}^{k'} \geq-\Delta^{t} \hbox{ if } i \in \mathcal{I}_P \cup \mathcal{I}_Q, i_F^{k'} \in \mathcal{I}_S,\\
-\Delta^{t} \text{ and } \bar{a}_{i_F^{k'}}^{k'} \geq -\Delta^{c}\hbox{ if } i\in \mathcal{I}_S, i_F^{k'} \in \mathcal{I}_P \cup \mathcal{I}_Q,\\
-\Delta^{t} \text{ and } \bar{a}_{i_F^{k'}}^{k'} \geq -\Delta^{t}\hbox{ if } i, i_F^{k'} \in \mathcal{I}_S.
\end{cases}
\end{align}
Note that the incentive condition ensures that a vehicle \( i \in \mathcal{I} \) prefers to change lanes if its expected acceleration on the target lane is significantly greater than its current acceleration. The safety condition guarantees that the lane-changing maneuver of vehicle \( i \in \mathcal{I} \) does not cause the expected following vehicle on the target lane and itself to brake excessively.\\
Now we introduce a lane-changing probability function.  The well-posedness of the Bando-FtL model has been recently proved in \cite{GongKeimer2022}. In particular, there exists $M \in \mathbb{R}_{\geq 0}$ such that for every $t \in [0, T]$ and $i \in \mathcal{I}$, the acceleration is bounded by M, i.e. $|a_i(t)| < M$. Let $\Delta = \min \{M, \Delta^{\gamma}\}_{\gamma \in \{c, t, cc, tc, ct, tt\}}$. Define three probability distributions as follows, 
\begin{align}
\label{eqn: probability_function}
  &p_j:(\mathbb{R}_{\geq 0})^5\rightarrow [0,1];\\ \nonumber 
  &p_j (b_1, b_2, b_3, b_4, b_5)= \frac{1}{C_j}\big(1- e^{- \gamma_j b_1b_2b_3b_4b_5}\big),
  \end{align}
where $\gamma_j>0$, $C_j$ are renormalization constants given by 
$$C_j=\max_{[0, 2M-\Delta]^5}\big(1- e^{-\gamma_j b_1b_2b_3b_4b_5}\big) = 1-e^{-\gamma_j (2M-\Delta)^5},$$ and $j=1, 2, 3$. The third probability function $p_3$ will be used in the infinite-dimensional case.  

A car or an autonomous vehicle $i\in \mathcal{I}_P\cup \mathcal{I}_Q$ will perform lane-change from lane $k\in \mathcal{K}$ to lane $k' \in \{k-1, k+1\} \cap \mathcal{K}$ under both the incentive and safety conditions with a probability given by
\begin{align*}
p_1\Big([\bar{a}_i^{k'}-a_i^k-\Delta^{cc}]_{+},\,[\bar{a}_i^{k'}-a_i^k-\Delta^{tc}]_{+},
\,[\bar{a}_i^{k'}+\Delta^{c}]_{+},\,[\bar{a}_{i_F^{k'}}^{k'}+\Delta^{c}]_+,[\bar{a}_{i_F^{k'}}^{k'} +\Delta^{t}]_{+}\Big).
\end{align*}
The probability of a truck $i\in \mathcal{I}_S$ performing lane-change is 
\begin{align*}
p_2\Big([\bar{a}_i^{k'}-a_i^k-\Delta^{ct}]_{+},\,[\bar{a}_i^{k'}-a_i^k-\Delta^{tt}]_{+},\,[\bar{a}_i^{k'}+\Delta^{t}]_{+},\,[\bar{a}_{i_F^{k'}}^{k'}+\Delta^{c}]_+,[\bar{a}_{i_F^{k'}}^{k'} +\Delta^{t}]_{+}\Big).
\end{align*} 

Here, \([r]_{+}\) represents the positive part of \(r \in \mathbb{R}\). That is, 
\[
[r]_{+} = 
\begin{cases} 
r & \text{if } r > 0, \\
0 & \text{if } r \leq 0.
\end{cases}
\]
 
The probability functions suggest that a vehicle is more likely to initiate a lane change when the expected acceleration in the target lane substantially exceeds its current acceleration, and when the maneuver involves less excessive braking for both the vehicle and its new follower. Observe that for modeling purposes, it is not relevant to assign a rule specifying if a vehicle is moving to the lane on the left or on the right, but only if a lane-change may happen or not.\\

\textbf{ The infinite-dimensional case}\\
A scenario involving an equal number of autonomous vehicles, cars and trucks is futuristic, however recent experiments aim to influence traffic with a large amount of human-driven vehicles by introducing few autonomous vehicles \cite{Kardous}. Therefore, we consider $Q$ autonomous vehicles and infinitely many cars and trucks.  The behaviour of the single human-driven vehicle does not play anymore a relevant role, what makes the difference is the average behaviour of the cars and trucks, for this reason we investigate the lane-changing behavior of all vehicles by looking into the following average acceleration of cars $A_P^k$, and the average acceleration of trucks $A_S^k$, on lane $k \in \mathcal{K}$: 
\begin{align*}
&A_P^k= H_1^{cc}\ast_1 (\mu^k_P+\mu^k_Q)+H_1^{tc}\ast_1\mu^k_S+H_2^{cc}\ast(\mu^k_P+\mu^k_Q)+H_2^{tc}\ast\mu^k_S,\\
&A_S^k =H_1^{ct}\ast_1 (\mu^k_P+\mu^k_Q)+H_1^{tt}\ast_1\mu_S^k+H_2^{ct}\ast(\mu^k_P+\mu^k_Q)+H_2^{tt}\ast\mu_S^k. 
\end{align*}
The probability of cars changing from lane $k \in \mathcal{K}$ to lane $k'\in \{k-1, k+1\}\cap \mathcal{K}$ lane is 
\begin{align*}
p_1\Big([A_P^{k'}-A_P^k-\Delta^{cc}]_{+},\,[A_S^{k'}-A_P^k-\Delta^{tc}]_{+},\,[A_P^{k'}+\Delta^{c}]_{+},
\,[A_S^{k'}+\Delta^{t}]_{+}, \,[A_S^{k'}+\Delta^{t}]_{+}\Big),
\end{align*}
the probability of trucks performing lane-change from lane $k \in \mathcal{K}$ to lane $k'\in \{k-1, k+1\}\cap \mathcal{K}$ is 
\begin{align*}
p_2\Big([A_P^{k'}-A_S^k-\Delta^{ct}]_{+},\,[A_S^{k'}-A_S^k-\Delta^{tt}]_{+},\,[A_P^{k'}+\Delta^{c}]_{+},
\,[A_S^{k'}+\Delta^{t}]_{+}, \,[A_S^{k'}+\Delta^{t}]_{+}\Big),
\end{align*}
and finally the probability of an autonomous vehicle $i \in \mathcal{I}_Q$ performing lane-change from lane $k \in \mathcal{K}$ to lane $k'\in \{k-1, k+1\}\cap \mathcal{K}$ is 
\begin{align*}
p_3\Big([\bar{a}_i^{k'}-A_P^k-\Delta^{cc}]_{+},\,[\bar{a}_i^{k'}-A_S^k-\Delta^{tc}]_{+},\,[\bar{a}_i^{k'}+\Delta^{c}]_{+},\,[A_S^{k'}+\Delta^{t}]_{+}, \,[A_P^{k'}+\Delta^{t}]_{+}\Big),
\end{align*}
where the probability functions $p_j$, $j=1, 2, 3$ are defined in \cref{eqn: probability_function}.
\subsection{Cool-down time}
\label{sec:cool_down_time}
Now we introduce the model assumption,  {cool-down} time as in \cite{gpv_meanfield}, which is critical to describe the frequencies of the vehicles' lane-changing behavior and to prove the well-posedness of our heterogeneus multi-lane traffic model.

By empirical observations, the lane-changing frequency of vehicles on the highway is low. For instance, a study analyzing a two dimensional dataset recorded on a German highway shows that only $15\%$ of the vehicles performed lane-change while traveling the recorded road segment  (see \cite{HertyVisconti2018}). For this reason, the chance of two vehicles performing lane-change at exactly the same time is even lower. Therefore, it is reasonable to assume that there are not two vehicles changing lane at the same time on the entire road. To achieve this, we associate each vehicle $i \in \mathcal{I}$ a timer $\tau_i$ and assume that the initial timers for two different vehicles are different. We also introduce the  {cool-down} time $\bar{\tau}=\tfrac{T}{N_{\tau}}$, where $N_{\tau} \in \mathbb{N}_{\geq 0}$ is large and assume that vehicle $i\in \mathcal{I}$ checks the lane-changing conditions only when its timer reaches the cool-down time, $\bar{\tau}$. 

In addition, we reset the vehicle's timer to $0$ once its timer reaches the cool-down time $\bar{\tau}$. Specifically, for each vehicle $i \in \mathcal{I}$, its timer $\tau_i$ satisfies the following 
\[  \dot{\tau}_i(t)=1, \, \tau_i(0)=\tau_{i,0}, \, t \in [0, \bar{\tau})\]
where $\tau_{i_1, 0} \not =\tau_{i_2, 0}$ if $i_1 \not = i_2 \in \mathcal{I}$. Note that one can also model large lane-changing frequencies by choosing small cool-down time ~$\bar{\tau}$. In the case of finitely many vehicles, the presence of the cool-down time, $\bar{\tau}$ allows us to consider a small time interval $[0, t_1]$ when there is no vehicle changing lane. Similarly, in the case of infinitely many cars and trucks but finitely many autonomous vehicles, due to the definition of the cool-down time, there is a small time interval $[0, t_2]$ when there is no autonomous vehicle changing lane. In particular, \(t_1 = \min_{i \in \mathcal{I}} \{\bar{\tau} - \tau_{i, 0}\} \text{ and } t_2 = \min_{i \in \mathcal{I}_{Q}} \{\bar{\tau} - \tau_{i, 0}\}.\)

{
\subsection{Generalization to multi-population models }
\label{sub_sec_multi_pop}
In this subsection, we derive a general formulation for both finite-dimensional and infinite-dimensional hybrid systems in the case of $M+1$ vehicle populations, comprising $M$ types of human-driven vehicles (such as cars, trucks, buses, motorcycles, and subclasses thereof) and autonomous vehicles. We briefly illustrate this case, as the derivation closely resembles the one developed for the simpler setting with only cars, trucks, and AVs. Indeed, the theory developed in the preliminary works of Gong, Piccoli, and Visconti \cite{gpv_meanfield, 9302751} can be naturally extended to more complex scenarios by applying the same techniques.

 Consider $A_k$ the number of autonomous vehicles on the lane $k\in\lbrace 1,\dots, L\rbrace$, with $\mathcal{I}_0^k$ the correspondent sets of indices. On the other hand, let $V_{n,k}$ with $n=1,\dots,M$ be the number of human-driven vehicles in the lane $k$ with indices in $\mathcal{I}_n^k$. 
The dynamic of this multi-population frame is given by the following system of first order ODE representing the Bando-FtL model in convolutional form:\vspace{-0.2cm}
 \begin{align}
 \label{gen_fin_dim_sys}
     &\dot{x}_i=v_i, \quad i \in \mathcal{I},\notag\\
  &\dot{v}_i=\begin{cases} \Big(\scalebox{0.9}{$\sum\limits_{m=0}^M H^{mn}_1*_1 \mu_m^k+\sum\limits_{m=0}^M H^{mn}_2* \mu_m^k$}\Big)(x_i, v_i)+u_i
      \quad i \in \mathcal{I}_0^k,
     \\
     \\
     \Big(\scalebox{0.9}{$\sum\limits_{m=0}^M H^{mn}_1*_1 \mu_m^k+\sum\limits_{m=0}^M H^{mn}_2* \mu_m^k$}\Big)(x_i, v_i)
     \quad i \in \mathcal{I}_n^k.
     \end{cases}
    \end{align}
 Here $\mu_n^k$, with $n\in\lbrace 0,\dots, M\rbrace$, represents again the atomic measures supported on absolutely continuous trajectories $t \in [0, T] \to (x_i(t), v_i(t)) \in \mathbb{R} \times \mathbb{R}_{\geq 0}$ for each type of vehicle in the lane $k$. Explicitly we have:\vspace{-0.2cm}
 \begin{align}
     \mu_0^k(t)&=\frac{1}{A_k(t)}\sum_{i\in\mathcal{I}_0^k(t)} \delta_{(x_i(t),v_i(t))},\notag\\ 
     \mu_1^k(t)&=\frac{1}{V_{1,k}(t)}\sum_{i\in\mathcal{I}_1^k(t)} \delta_{(x_i(t),v_i(t))},\notag\\
     \vdots\notag\\
     \mu_M^k(t)&=\frac{1}{V_{M,k}(t)}\sum\limits_{i \in \mathcal{I}_M^k(t)} \delta_{(x_i(t),v_i(t))}\notag.
 \end{align}
 The convolutional kernels in (\ref{gen_fin_dim_sys}) are defined in such a way that they can represents all the possible combinations of vehicles depending on their order, indeed 
 \vspace{-0.2cm}
 \[H_q^{mn}:\, \mathbb{R}\times\mathbb{R}_{\geq 0 }\rightarrow \mathbb{R} \quad \hbox{ with } q\in\lbrace 1,2 \rbrace\,, \quad m,n\in\lbrace 0,\dots, M\rbrace,\]
 and the structure of these maps is the same described for the two populations case in \ref{SecA}. The control appears only in the acceleration of autonomous vehicles and aims to influence the dynamic of the other populations. 
 Observe that in the case of $M=2$, i.e. with autonomous vehicles and two population of human-driven vehicles (cars and trucks), if we assume that the convolutional kernel for the first two types is the same, then we find again the dynamic \eqref{finite_dimensional_system}. 

 We point out that the lane changing conditions (Incentive and Safety), the probability for a vehicle of performing lane change and the definition of controlled hybrid system can be formulated in analogous way up to heavier notation. It can be also proved rigorously that the mean field limit (i.e, assuming that the number of vehiches for each class, apart AV, goes to $+\infty$) for (\ref{gen_fin_dim_sys}) is given by,
 \begin{align*}
 &\dot{x}_i=v_i, \quad i \in \mathcal{I}_0,\notag\\
 & \dot{v}_i =  \Big(\sum\limits_{m=0}^M H^{m0}_1*_1 \nu_m^k+\sum\limits_{m=0}^M H^{m0}_2* \nu_m^k\Big)(x_i, v_i)+u_i,
      \qquad i \in \mathcal{I}_0^k,\notag\\
      &\partial_t\nu^k_n+v\partial_x\nu^k_n+
      \partial_v\Big[\Big(\sum\limits_{m=0}^M H^{mn}_1*_1 \nu_m^k+\sum\limits_{m=0}^M H^{mn}_2* \nu_m^k\Big)\nu_n^k\Big]\\
     &\qquad \qquad \qquad= G_n(\nu_n^k, \nu_m^k, \nu_n^{k'}, \nu_m^{k'}),\qquad n\in \lbrace 1,\dots, M, \rbrace\notag
     \end{align*}
 with $\nu_n^k\in \mathcal{M}^{+}(\mathbb{R} \times \mathbb{R}_{\geq 0})$ representing the density  distribution of vehicles of class $n$ in the space of position and velocity.
  In the new source term, the input $\nu_m^k$ is the density of vehicles of the class $m\neq n$ entering the contiguous lane, therefore, it stands for $M-1$ densities (same for $\nu_m^{k'}$ ).
  
  This extension of the model to multiple populations allows for a fine analysis of traffic, including for example, small, medium, and large vehicles, motorcycles, bicycles and so on.}
\section{The optimal control problem on a finite-dimensional hybrid system}\label{opt_FD}
In this  section, we again consider the multi-lane and multi-population traffic with $P$ cars, $S$ trucks and $Q$ autonomous vehicles. The continuous dynamics of the finitely many vehicles without lane-change and the discrete events generated by the vehicles' lane-changing behaviors lead us to consider a finite-dimensional hybrid system $\Sigma_1$. We distinguish the autonomous vehicles from the others by adding control terms to their accelerations. 

Let $X=\mathbb{R} \times \mathbb{R}_{\geq 0} \times [0, \bar{\tau})$ and 
\begin{equation}\label{Ell}\mathcal{L} = \left\{\ell = (\ell_i)_{i\in\mathcal{I}} \in \mathcal{K}^{P+Q+S} \right\}\end{equation}
be the set of symbols that represent all possible lane labels of all vehicles including cars, trucks and autonomous vehicles. Additionally, before giving the definition of the finite-dimensional hybrid system $\Sigma_1$, we define the following two sets: the set $A_\ell$ containing the position-velocity-timer vectors of all vehicles among which there are at least two vehicles occupying the same lane and position at certain time, and the set $LC(\Sigma_1)$ representing the lane-changing mechanism of the finitely many vehicles: 
\begin{align}
\label{def_Al}
    A_{\ell} &= \Big\lbrace\big(x_{i}, v_{i}, \tau_{i}\big) _{i \in \mathcal{I}} \in X\colon \exists t\in [0, T], i_1, i_2 \in \mathcal{I}, s.t.,\notag\\
   \Big(& \scalebox{0.89}{$x_{i_1}(t)=x_{i_2}(t)\Big) \wedge \Big(\ell_{i_1}(t) = \ell_{i_2}(t)\Big), \text{ with } \ell_{i_1}, \ell_{i_2} \in \mathcal{K} \Big\rbrace,$} 
\end{align}\vspace{-3mm}
\begin{align}
&LC(\Sigma_1)=\Big\lbrace(\ell, (x_i, v_i, \tau_i), \ell', (x_i', v_i', \tau_i'))_{i\in \mathcal{I}} \in (\mathcal{L} \times X)^2\colon  \nonumber\\ 
& \qquad\qquad\exists \, i_0 \in \mathcal{I},\exists\, t_0 \in [0, \bar{\tau}), s.t.,
j \not = i_0, \Big((\ell_j(t_0), x_j(t_0), v_j(t_0),  \nonumber \\ 
&\qquad\qquad\tau_j(t_0))= (\ell_j'(t_0), x_j'(t_0), v_j'(t_0), \tau_j'(t_0))\Big) \wedge \Big(x_{i_0}(t_0), \nonumber \\
&\qquad\qquad v_{i_0}(t_0)) = (x_{i_0}'(t_0), v_{i_0}'(t_0))\Big), \ell'_{i_0}(t_0)=\ell_{i_0}(t_0)\pm 1, \nonumber \\
&\qquad\qquad\tau_{i_0}'(t_0)=0
\Big\rbrace. 
\label{def:switvhing_set}
\end{align}
Now we are ready to define the controlled finite-dimensional hybrid system. 
\begin{definition}[The controlled finite-dimensional hybrid system]
\label{def_finite_hybrid_system}
A finite dimensional hybrid system is a $6$-tuple $\Sigma_1 = (\mathcal{L}, \mathcal{M},U,\mathcal{U}, g, SW)$ where:
\vspace{0.5em}
\begin{enumerate}\setlength\itemsep{0.5em}
\item[(1)]  $ \mathcal{L}$ is a finite set of symbols representing all possible lane labels as defined in \eqref{Ell}. 
Here we call $\ell \in \mathcal{L}$ a location of the hybrid system $\Sigma_1$;

\item[(2)]  $\mathcal{N} =\{\mathcal{N}_{\ell}\}_{\ell \in \mathcal{L}}$, 
where $\mathcal{N}_{\ell} = (X \setminus A_{\ell})^{P+Q+S}$ is the space of position-velocity-timer vectors of all vehicles,
with $A_{\ell}$ defined as in \eqref{def_Al};
\item[(3)] $U=\{U_{\ell}\}_{\ell\in \mathcal{L}}$ represents the control space, $U_{\ell} =I^{Q}$, where $I \subset [-U_{\max}, U_{\max}]$ is compact with $U_{\max}>0$;

\item[(4)] $\mathcal{U} = \{\mathcal{U}_{\ell}\}_{\ell \in \mathcal{L}}$ is such that $\mathcal{U}_{\ell} = \{u: [0, T]\subset \mathbb{R}_0^+ \to U_{\ell\,} \mbox{ integrable}\}$  which represents the set of admissible controls at location $\ell$;
\item[(5)]  $g = \{g_{\ell}\}_{\ell \in \mathcal{L}}$ with
$g_{\ell} \colon \mathcal{N}_{\ell} \times \mathcal{U}_{\ell} \mapsto \mathbb{R}^{3(P+Q+S)}$, is such that for every $(x_i, v_i, \tau_i, u_i) \in \mathcal{N}_{\ell_i} \times \mathcal{U}_{\ell_i}$, it holds $g_{(\ell_i)}(x_i, v_i, \tau_i, u_i)=(v_i, a_i, 1)$, where
$a_i=\dot{v}_i$ is defined as in systems \eqref{finite_dimensional_system};
\item[(6)] $SW \text{ is a subset of } LC(\Sigma_1)$, where 
$LC(\Sigma_1)$ is the set of states for which a lane-changing can occur, that is (\ref{def:switvhing_set}).
\end{enumerate}
\end{definition}

For the well-posedness of the above finite-dimensional hybrid system, we refer to \cite{CGPmultipopulation}. Specifically, if the convolution kernels \( H_1^{\gamma}, H_2^{\gamma} \colon \mathbb{R} \times \mathbb{R}_{\geq 0} \to \mathbb{R} \), for \( \gamma \in \{cc, tc, ct, tt\} \), are locally Lipschitz and exhibit sub-linear growth, then the trajectories of all vehicles—including cars, trucks, and autonomous vehicles—are Lipschitz continuous in time over any interval during which no lane changes occur.

In this setting, we can introduce an optimal control problem associated with the finite-dimensional hybrid system of \cref{def_finite_hybrid_system} on the time interval $[0, t_1)$  where there is no vehicle changing lanes, which can be extended to the whole time interval $[0, T]$.\\ 
To state and prove the next results, we need to consider the product space 
\begin{equation*}
 \mathcal{X}\colon\,=(\mathbb{R} \times \mathbb{R}_{\geq 0})^M \times \mathcal{M}^{+}(\mathbb{R}\times \mathbb{R}_{\geq 0})^2  
\end{equation*}
 endowed with the following metric: for any $(x_1, v_1, \mu_1, \nu_1),(x_2, v_2, \mu_2,\nu_2)\in \mathcal{X}$,
\begin{align}\label{the metric_generalied}
&\|(x_1, v_1, \mu_1, \nu_1) - (x_2, v_2, \mu_2,\nu_2)\|_{\mathcal{X}} =\\
=&\sum \limits_{j=1}^{M}\scalebox{0.85}{ $(|x_{1, j}- x_{2, j}|+|v_{1, j}-v_{2,j}|)+\mathcal{W}_1^{1,1}(\mu_1, \mu_2)+\mathcal{W}_1^{1,1}(\nu_1, \nu_2)$},\notag
\end{align}
where $M \in \mathbb{N}_{0}$, $\mathcal{W}_1^{1,1}$ is the generalization of the standard Wasserstein distance introduced in \cite{piccoli2014generalized}.
\begin{definition}[Optimal control problem associated with a finite-dimensional hybrid system]
\label{fin_dim_functional}
Find $u^{*} \in L^1([0, t_1); I)^Q$, such that\vspace{-3mm}
\begin{equation}\label{eqn_FPSk}
 F_{P, S}(u^{*}) = \min\limits_{u\in L^1([0, {t_1}), I)^Q} F_{P, S}(u).
    \end{equation}
    \vspace{-1mm}
    with the functional $ F_{P, S}$ given by
     \vspace{-1mm}
   \begin{equation}
\label{eqn_FPS}
\begin{aligned}
     \sum\limits_{k \in \mathcal{K}} \int_0^{t_1} \scalebox{0.88}{$\left\{L_k(x^k(t), v^k(t), \mu_{P}^k(t),\mu_{S}^k(t))+\frac{1}{Q_k(t)} \sum \limits_{j=1}^{Q_k(t)} |u_{j}^k(t)|\right\}\mathrm{d} t$}
\end{aligned}
\end{equation}
where $(x^k, v^k, \mu_P^k, \mu_S^k) \in (\mathbb{R} \times \mathbb{R}_{\geq 0})^{Q_k(t)} \times \mathcal{P}(\mathbb{R} \times \mathbb{R}_{\geq 0})^2$ are solutions to the  finite dimensional hybrid system \eqref{def_finite_hybrid_system} on the time interval $[0, t_1)$ and $L_k:(\mathbb{R} \times \mathbb{R}_{\geq 0})^{Q_k(t)} \times \mathcal{P}(\mathbb{R} \times \mathbb{R}_{\geq 0})^2\to \mathbb{R}$ are running costs continuous with respect to the metric \eqref{the metric_generalied}.
\end{definition}
Notice that the cost functional defined in (\ref{eqn_FPS}) contains two pieces of different nature: the former is a classic Lagrangian involving position and velocity of autonomous vehicles and the density measure of cars and trucks; the latter is a weighted $L^{1}$-norm of the vector of controls. 
For \eqref{fin_dim_functional} we can prove the following result:

\begin{theo}\label{OgammaD}
    The finite horizon optimal control problem \eqref{eqn_FPSk}  associated to the finite dimensional hybrid system has solutions.
\end{theo}
\begin{proof}
    Let $(u_n)_{n\in\mathbb{N}}$ be a minimizing sequence realizing at its limit the minimum of the cost functional \eqref{eqn_FPS}. Such a sequence must be bounded in $L^1([0, t_1); I)^Q$, hence it admits a subsequence that we keep indicating $(u_n)_{n\in\mathbb{N}}$ which is weakly convergent to a $u^*\in L^1([0, t_1); I)^Q$. At the same time, for a fixed control $u^n$, the solutions $(x^n_i (t),v^n_i(t))$ with $i\in\mathcal{I}_P^k\cup\mathcal{I}_Q^k\cup\mathcal{I}_S^k$ of  \eqref{finite_dimensional_system} are equibounded and equi-Lipschitz continuous in time, therefore they converge to trajectories $(x^*_i (t),v^*_i(t))$ which are solutions to \eqref{finite_dimensional_system} for the control $u^*$. Let 
    \begin{align}
    \mu_P^{n,k}(t) &= \frac{1}{P_k(t)} \sum_{i \in \mathcal{I}_P^k(t)} \delta_{(x^n_i(t), v^n_i(t))},\qquad \qquad \mu_P^{*,k}(t) &= \frac{1}{P_k(t)} \sum_{i \in \mathcal{I}_P^k(t)} \delta_{(x^*_i(t), v^*_i(t))}, \\
    \mu_S^{n,k}(t) &= \frac{1}{S_k(t)} \sum_{i \in \mathcal{I}_S^k(t)} \delta_{(x^n_i(t), v^n_i(t))},\qquad \qquad \mu_S^{*,k}(t) &= \frac{1}{S_k(t)} \sum_{i \in \mathcal{I}_S^k(t)} \delta_{(x^*_i(t), v^*_i(t))},
\end{align}
and denote by $({\bf{x}}^{n,k}(t),{\bf{v}}^{n,k}(t))\in(\mathbb{R} \times \mathbb{R}_{\geq 0})^{Q_k(t)}$ the vector of trajectories of AVs on lane $k$ associated to the control $u^n$ and by $({\bf{x}}^{*,k}(t),{\bf{v}}^{*,k}(t))\in(\mathbb{R} \times \mathbb{R}_{\geq 0})^{Q_k(t)}$ the vector of trajectories of AVs on lane $k$ for the limit control $u^*$. The uniform convergence of the trajectories and their compact support also allow us to conclude that 
$$\lim_{n\to\infty}({\bf{x}}^{n,k}(t),{\bf{v}}^{n,k}(t),\mu_P^{n,k}(t),\mu_S^{n,k}(t) )=({\bf{x}}^{*,k}(t),{\bf{v}}^{*,k}(t),\mu_P^{*,k}(t),\mu_S^{*,k}(t) )$$
hence,
     \begin{equation}
         \lim\limits_{n\to \infty}\int_0^{t_1}L_k({\bf{x}}^{n,k}(t),{\bf{v}}^{n,k}(t),\mu_P^{n,k}(t),\mu_S^{n,k}(t) )dt =\int_0^{t_1}L_k({\bf{x}}^{*,k}(t),{\bf{v}}^{*,k}(t),\mu_P^{*,k}(t),\mu_S^{*,k}(t) )dt.
    \end{equation}

The weak convergence of $(u_n)_{n\in\mathbb{N}}$ to $u^*\in L^1([0, t_1); I)^Q$ implies the lower semicontinuity of the weighted $L^{1}$-norm of the vector of controls, i.e.
$$\liminf_{n\to\infty}\frac{1}{Q_k(t)} \sum \limits_{j=1}^{Q_k(t)} |u_{j}^{n,k}(t)|dt\geq \frac{1}{Q_k(t)} \sum \limits_{j=1}^{Q_k(t)} |u_{j}^{*,k}(t)|dt.$$

We can conclude that $u^*$ is an optimal control by the direct method of calculus of variations.
\end{proof}

\section{The optimal control problem on an infinite-dimensional hybrid systems}\label{inf_dim}
In the following, we introduce a definition of controlled infinite-dimensional hybrid system. The rigorous proof of the passage from the finite to the infinite system is given in \cite{CGPmultipopulation}. Here again the system describes the dynamics of multi-lane traffic with cars, trucks, and autonomous vehicles. The difference is that the evolution for the first two types of vehicle is given in terms of their density, while the number of autonomous vehicles is still finite. The lane changing behavior of the autonomous vehicles generates discrete events of the infinite dimensional hybrid system. Let $(x_i, v_i)$ be the position-velocity vector of the autonomous vehicle $i \in \mathcal{I}_Q$ and $\nu_c^k, \nu_t^k \in \mathcal{M}^{+}(\mathbb{R} \times \mathbb{R}_{\geq 0})$
the density distribution of cars 
and trucks on the $k$-th lane in space and velocity with $\mu_Q^k$ 
the empirical measure for autonomous vehicles defined in (\ref{emAV}). The continuous dynamics of all vehicles are given by
\begin{align}\label{infinite_dimensional_system}
&\dot{x}_i=v_i \quad i \in \mathcal{I}_Q\notag\\
& \dot{v}_i =  \left(H_1^{cc}\ast_1 (\nu_c^k+\mu_Q^k)+H_1^{tc}\ast_1\nu_t^k\right)(x_i,v_i)\notag\\
  &\qquad+\left(H_2^{cc}\ast(\nu_c^k+\nu_Q^k)+H_2^{tc}\ast\nu_t^k\right)(x_i,v_i) +u_i
     \quad i \in \mathcal{I}_Q^k\notag\\
    &\partial_t\nu^k_c+v\partial_x\nu^k_c+
     \partial_v\Big[\big(H_1^{cc}\ast_1 (\nu_c^k+\nu_Q^k)+H_1^{ct}\ast_1\nu^k_t\notag\\
     &\qquad +H_2^{cc}\ast(\nu_c^k+\nu_Q^k)+H_2^{ct}\ast\nu^k_t\big)\nu_c^k\Big]=G_1(\nu_c^k, \nu_t^k, \nu_c^{k'}, \nu_t^{k'})\notag\\
  &\partial_t\nu^k_t+v\partial_x\nu^k_t+\partial_v\Big[\big(H_1^{ct}\ast_1 (\nu_c^k+\nu_Q^k)+H_1^{tt}\ast_1\nu^k_t\notag\\& \qquad +H_2^{ct}\ast(\nu_c^k+\nu_Q^k)+H_2^{tt}\ast\nu^k_t\big)\nu^k_t\Big]=G_2(\nu_c^k, \nu_t^k, \nu_c^{k'}, \nu_t^{k'}),
     \end{align}
where $k \in \mathcal{K}$, the source terms $G_1$ and $G_2$ are generated by the lane-changing behavior of the cars and trucks. To define the infinite-dimensional controlled hybrid system  we consider $X=\mathbb{R} \times \mathbb{R}_{\geq 0} \times [0, \bar{\tau})$ and $\Tilde{\mathcal{L}} = \left\{\ell = (\ell_i)_{i\in\mathcal{I}_Q} \in \mathcal{K}^{Q} \right\}$ the set of symbols representing all the possible lane labels for autonomous vehicles. Moreover we introduce $\Tilde{A}_{\ell}$ set of triples position-velocity-timer  and $LC(\Sigma_2)$ representing the lane-changing mechanism of the finitely many autonomous vehicles:
\begin{align}
\label{def_Ali}
    \Tilde{A}_{\ell} &= \Big\lbrace\scalebox{0.89}{$\big(x_{i}, v_{i}, \tau_{i}\big) _{i \in \mathcal{I}_Q}\in X\colon \exists t\in [0, T], i_1, i_2 \in \mathcal{I}_Q,$} \notag\\
    &\scalebox{0.89}{$s.t., \Big(x_{i_1}(t)=x_{i_2}(t)\Big)\wedge \Big(\ell_{i_1}(t) = \ell_{i_2}(t)\Big), \hbox{ with }\ell_{i_1},\ell_{i_2}\in \mathcal{K} \Big\rbrace,$}\\\notag
\label{def:switvhing_set2}
&\scalebox{0.95}{$LC(\Sigma_2)=
\Big\lbrace\big(\ell, (x_i, v_i, \tau_i), \ell', (x_i', v_i', \tau_i')\big)_{i\in \mathcal{I}_Q}\in (\Tilde{\mathcal{L}} \times X)^2\colon$} \notag \\ 
\exists& \scalebox{0.90}{$\, i_0 \in \mathcal{I}_Q,\exists\, t_0 \in [0, \bar{\tau}), s.t.,
j \not = i_0, \Big((\ell_j(t_0), x_j(t_0), v_j(t_0), \tau_j(t_0))$} \notag\\ 
=& (\ell_j'(t_0), x_j'(t_0), v_j'(t_0), \tau_j'(t_0))\Big) \wedge \Big((x_{i_0}(t_0), v_{i_0}(t_0))\notag \\
=& (x_{i_0}'(t_0), v_{i_0}'(t_0))\Big), \ell'_{i_0}(t_0)=\ell_{i_0}(t_0)\pm 1, \tau_{i_0}'(t_0)=0
\Big\rbrace. 
\end{align}
\begin{definition}[The controlled infinite-dimensional hybrid system]
\label{def_infinite_hybrid_system}
A infinite dimensional hybrid system is a $6$-tuple  $\Sigma_2 = ( \Tilde{\mathcal{L}}, \mathcal{M},U,\mathcal{U}, g, SW)$
\vspace{0.5em}
\begin{enumerate}\setlength\itemsep{0.5em}
\item[(1)]  $ \Tilde{\mathcal{L}} = \left\{\ell = (\ell_i)_{i\in\mathcal{I}_Q} \in \mathcal{K}^{Q} \right\}$ is a finite set of symbols that represent all possible lane labels of the autonomous vehicles;
\item[(2)]  $\mathcal{N} =\{\mathcal{N}_{\ell}\}_{\ell \in \Tilde{\mathcal{L}}}$, 
where $\mathcal{N}_{\ell} = (X \setminus \Tilde{A}_{\ell})^{Q}\times(\mathcal{M}^+(\mathbb{R}\times\mathbb{R}_{\geq 0}))^{2L}$,
with $\Tilde{A}_{\ell}$ given by \eqref{def_Ali};
\item[(3)] $U=\{U_{\ell}\}_{\ell\in \Tilde{\mathcal{L}}}$ represents the control space; 

\item[(4)] $\mathcal{U} = \{\mathcal{U}_{\ell}\}_{\ell \in \Tilde{\mathcal{L}}}$ is such that $\mathcal{U}_{\ell} = \{u: \text{Dom}(u)\subset \mathbb{R}_0^+ \to U_{\ell} \mbox{   integrable}\}$  which represents the set of admissible controls at location $\ell$;
\item[(5)]  $g = \{g_{\ell}\}_{\ell \in \tilde{\mathcal{L}}}$,
$g_{\ell} \colon \mathcal{N}_{\ell} \times \mathcal{U}_{\ell} \mapsto \mathbb{R}^{3Q}$, such that for every $(x_i,v_i, \tau_i, \nu_c, \nu_t, u_i) \in \mathcal{N}_{\ell_i} \times \mathcal{U}_{\ell_i}$, $g_{(\ell_i)}(x_i,v_i, \tau_i, \nu_c, \nu_t, u_i)=(v_i, a_i, 1)$, where
$a_i=\dot{v}_i$, $i \in \mathcal{I}_Q$ as defined in systems \eqref{infinite_dimensional_system};
\item[(6)] $S \text{ is a subset of } LC(\Sigma_2)$, where 
$LC(\Sigma_2)$ is the set of states for which a lane-changing can occur, that is (\ref{def:switvhing_set2}).
\end{enumerate}
\end{definition}
The well-posedness of the infinite-dimensional hybrid system is established in \cite{CGPmultipopulation, gpv_meanfield}.

\begin{definition}[Optimal control problem associated with an infinite-dimensional hybrid system]\label{Opt_contr_inf} Find $u^{*} \in L^1([0, t_2); I)^Q$, such that \vspace{-1mm}
\begin{equation}\label{infin_dim_functional}
 F(u^{*}) = \min\limits_{u\in L^1([0, t_2), I)^Q} F(u),
    \end{equation}
    and the functional $F$ is given by
    \vspace{-1mm}
\begin{equation}
\label{eqn_Fk}
\begin{aligned}
   \sum\limits_{k \in J} \int_0^{t_2} \scalebox{0.9}{$\left\{L_k(x^k(t), v^k(t), \nu_{c}^k(t),\nu_{t}^k(t))+\frac{1}{Q_k(t)} \sum \limits_{j=1}^{Q_k(t)} |u_{j}^k(t)|\right\}\mathrm{d} t$}
\end{aligned}
\end{equation}

where $(x^k, v^k, \nu_c^k, \nu_t^k)\in (\mathbb{R} \times \mathbb{R}_{\geq 0})^{Q_k(t)} \times \mathcal{M}^+(\mathbb{R} \times \mathbb{R}_{\geq 0})^2$ are solutions to system \eqref{infinite_dimensional_system} and each function $L_k:(\mathbb{R} \times \mathbb{R}_{\geq 0})^{Q_k(t)} \times \mathcal{P}(\mathbb{R} \times \mathbb{R}_{\geq 0})^2\rightarrow \mathbb{R}_{\geq 0}$ is continuous with respect to the metric defined in \eqref{the metric_generalied}.
\end{definition}

\begin{definition}
Let $X$ be a separable metric space and consider the sequence of functionals $F_N \colon X \mapsto (-\infty, \infty]$, $N\in \mathbb{N}$. Then $F_N$ $\Gamma$- converges to $F \colon X \mapsto (-\infty, \infty]$ if the following conditions are satisfied: 
\begin{itemize}
\item[1.](\,{\bf Lim inf inequality}\,)   For every  $u \in X$  and every sequence $ u_N \to u$, $$F(u) \leq \liminf \limits_{N \to \infty} F_N(u_N);$$
\item[2.](\,{\bf Lim sup inequality}\,) For every $ u \in X$,  there exists a sequence $ u_N \to u$, s.t. 
$$ F(u) \geq \limsup\limits_{N \to \infty} F_N(u_N).$$ 
\end{itemize}
\end{definition}

\begin{theo}
\label{thm_gamma}
For every lane $k \in \mathcal{K}$, the sequence of functionals $(F_{P,S}^k)_{P, S \in \mathbb{N}^{+}}$ on $L^1([0, t_2), \mathcal{U})^Q$ as defined in \eqref{eqn_FPS} $\Gamma$-converges to the functional $F^k$  in \eqref{eqn_Fk}.
\end{theo}
\begin{proof}
We start showing that the lim inf inequality is satisfied. Let $(u_N)_{N\in\mathbb{N}}$ be a sequence in {$L^1([0, t_2))^Q$} such that $u_{N} \rightharpoonup u \in L^1([0, t_2))^Q$. Furthermore, we require that as $N \to \infty$, the number of human-driven vehicles for each type goes to infinity as well (i.e. $P, S \to \infty$).
{ We claim that} for each $k \in \mathcal{K}$, 
\begin{equation}\label{LMP} F^k(u) \leq \liminf \limits_{P \to \infty,\, S \to \infty} F_{P,S}^k(u_N).\end{equation}
For each $N \in \mathbb{N}$, there exists a unique solution $({\bf{x}}^{N,k}(t),{\bf{v}}^{N,k}(t),\mu_P^{N,k}(t),\mu_S^{N,k}(t) ) \in (\mathbb{R}\times \mathbb{R}_{\geq 0})^{Q_k(t)} \times \mathcal{M}^{+}(\mathbb{R} \times \mathbb{R}_{\geq 0})$, $k \in \mathcal{K}$, to the finite-dimensional hybrid system.  
The key property needed to prove (\ref{LMP}) is 
\begin{equation*}
    \lim \limits_{N \to \infty}({\bf{x}}^{N,k}(t),{\bf{v}}^{N,k}(t),\mu_P^{N,k}(t),\mu_S^{N,k}(t) ))=({\bf{x}}^{k}(t),{\bf{v}}^{k}(t), \nu_{c}^{k}(t),\nu_{t}^{k}(t)) \end{equation*}
where $({\bf{x}}^{k}(t),{\bf{v}}^{k}(t), \nu_{c}^{k}(t),\nu_{t}^{k}(t))$ is the solutions to \eqref{infinite_dimensional_system} with control $u$, the limit is in the metric defined in (\ref{the metric_generalied}), and it holds by construction. This indeed implies that 
\begin{align*}
        \lim\limits_{N\to \infty}&\int_0^{t_0}L_k({\bf{x}}^{N,k}(t),{\bf{v}}^{N,k}(t),\mu_P^{N,k}(t),\mu_S^{N,k}(t) ))dt \notag\\
        &=\int_0^{t_0}L_k({\bf{x}}^{k}(t),{\bf{v}}^{k}(t), \nu_{c}^{k}(t),\nu_{t}^{k}(t))dt
\end{align*}
which together the lower-continuity of the $L^1$-norm gives (\ref{LMP}).\\
To prove the Lim sup inequality, let us fix $u_N\equiv u^*$ for all $N\in\mathbb{N}$. As before, we can associate to this sequence of controls a sequence of solutions $({\bf{x}}^{N,k}(t),{\bf{v}}^{N,k}(t),\mu_P^{N,k}(t),\mu_S^{N,k}(t) )$ to the finite dimensional hybrid system \eqref{finite_dimensional_system} convergent to $({\bf{x}}^{*,k}(t),{\bf{v}}^{*,k}(t),\mu_P^{*,k}(t),\mu_S^{*,k}(t))$ in the metric defined in (\ref{the metric_generalied}). By continuity of $L_k$ and since the sequence $(u_N)_{N\in\mathbb{N}}$ is constant we have 
$$F^k(u^*) = \limsup \limits_{P \to \infty,\, S \to \infty} F_{P,S}^k(u_N)$$
\end{proof}

\begin{cor}
The optimal control problem in Definition (\ref{Opt_contr_inf}) has solutions. 
\end{cor}
\begin{proof}
By \eqref{OgammaD}, for each $N \in \mathbb{N}$, there is an optimal control $u_{*, N}$ for the finite dimensional system. Note that since the sequence of controls $(u_{*, N})$ is bounded in $L^1([0, t_2),I)^Q$ which is compact in the weak topology, there exists a subsequence (for simplicity, we still use the same { notation}), such that $u_{*, N} \rightharpoonup u^{*} \in L^1([0, t_2),I)^Q$. By Theorem 7.8 in \cite{DalMaso}, $u^{*} $ minimizes $F$.  
\end{proof}
{
\begin{rem}
    Our results are presented within a non-hybrid analytical framework. We recall that each lane change is subject to a cool-down time, which enforces a minimal dwell time between two consecutive switches of a vehicle. This assumption prevents chattering phenomena and guarantees that only finitely many discrete transitions can occur within any finite time interval. Consequently, the continuous part of the dynamics can be analyzed in a non-hybrid setting, and all existence results established in this work remain valid in the presence of lane changes under this condition.
\end{rem}
\begin{rem}
     Similarly to what is been done for cars, trucks and autonomous vehicles, we can associate an optimal control problem both to the finite and infinite dimensional hybrid system and study the well posedness for the model involving \(M+1\) types of vehicles introduced in subsection \ref{sub_sec_multi_pop}. 
\end{rem}
}
\section{Microscopic Hybrid Multi-Population Simulations}\label{sec:simulations}

We aim to investigate the impact of the truck penetration rate on overall velocity variation within a hybrid, multi-population traffic flow model. To this end, we develop a microscopic traffic simulator that incorporates multiple vehicle populations and discrete lane-changing behavior. In this section, we outline the general simulation settings and provide a schematic overview of the numerical approach.

\subsection{Model Equations}
The microscopic traffic model relies on a slight alteration of the Bando-FtL terms; see Equation \ref{eq:bando-ftl} where the optimal velocity function is given by Equation \ref{eq:opt-velocity}. In the main Bando-FtL terms, the subscript, $n$ indexes the vehicles on the road, where $n$ refers to the ego vehicle and $n+1$ refers to the leader of the ego vehicle. Additionally, $h_n$ is the gap headway for vehicle $n$, meaning it is the distance from the front bumper of the vehicle$_n$ to the rear bumper of vehicle$_{n+1}$. Within the optimal velocity function, we require several parameters: $h$ the current headway, $d_s$ the ``safe following distance'' for the road, and $l_v$ the length of the vehicle.
\begin{align}\label{eq:bando-ftl}
    \dot{x}_n &= v_n \\
    \dot{v}_n &= \alpha_n(V(h_n) - v_n) + \beta_n \frac{v_{n+1} - v_n}{(h_n)^2}
\end{align}
\begin{equation}\label{eq:opt-velocity}
    V(h) = v_{max} \frac{\tanh(h-d_s) + \tanh(l_v + d_s)}{1 + \tanh(l_v + d_s)}
\end{equation}

\subsection{Simulation Setting}
To examine the effect of different proportions of the multi-populations, we designed a simulator for a multi-lane ring road (equal lane length, wrapping effect) with multiple populations of vehicles. Each scenario simulation has three lanes of equal length and two vehicle populations, ``cars'' and ``trucks''; each vehicle is able to lane change according to the incentive and safety conditions outlined in Equation \ref{lane-changing_acc}. Cars are implemented as vehicles with a length of 4.5 meters, whereas trucks are implemented as vehicles with a length of 13.6 meters; additionally, trucks have slower acceleration (reduced $\alpha, \beta$ parameters in Equation \ref{eq:bando-ftl}). Model parameters for Bando-FtL are based on values from previous work as in \cite{blondin_feedback_2019}.

\subsection{Algorithmic Flow}
The simulator contains two main components: a lane-changing section and a fixed step-size Runge-Kutta integrator. The general computational flow is described in Algorithm \ref{alg:main-comp-flow}. The simulator is implemented using an object-oriented programming (OOP) paradigm, with objects for vehicles and the road, where different populations of vehicles are subclasses of a Vehicle superclass. Additionally, the Road objects contain lanes, which are lists of vehicles -- each lane is implemented as a doubly-linked list for straightforward access to leader and follower vehicles for each ego vehicle. To compute headways for the computation of the Bando-FtL term (Equation \ref{eq:bando-ftl}), we iterate through the linked lists to obtain each vehicle's leader. Lane-changing behavior relies on the incentive and safety conditions for lane-changing as in Expressions \ref{exp:lane-changing}. The intuitive explanation for the incentive condition is the following: if the ego vehicle's acceleration in the new lane would be greater than its current acceleration plus some offset, then the ego vehicle is \textit{incentivized} to change lanes. The safety condition is similar, except we ensure that the acceleration within the new lane is above some minimum safety threshold; additionally, we ensure that the new follower of the ego vehicle does not exceed some maximum deceleration threshold.
\begin{align}\label{exp:lane-changing}
    \text{Incentive}:& \;\; \bar{a}_i^{k'}(t) \geq \bar{a}_i^k(t) + \Delta\\
    \text{Safety}:& \;\; \bar{a}_i^{k'}(t) \geq -\Delta \text{  and  } \bar{a}_{i_F^{k'}}^{k'}(t) \geq -\Delta
\end{align}
\begin{algorithm}
    \caption{Main computational flow of the simulator}\label{alg:main-comp-flow}
    \begin{algorithmic}
        \State $t \gets 0$
        \While {$t < t_{end}$} 
            \For {each vehicle, $veh$}
                \State lc $\gets$ Check-for-lane-change($veh$)
                \If{lc is True}
                    \State $veh \gets $ update-vehicle-info()
                    \State $lane \gets $ update-lane-info()
                \EndIf
            \EndFor
            
            \State Runge-Kutta-Update()
            \State $t \gets t + t_{step}$
        \EndWhile
    \end{algorithmic}
\end{algorithm}

\begin{algorithm}
    \caption{The computational flow for the lane-changing function}\label{alg:lane-change-comp-flow}
    \begin{algorithmic}
        \Procedure{Check-for-lane-change}{$veh$}
            \State successfulLaneChange $\gets$ False
            \State cooldown $\gets$ veh.getTimeOfLastLaneChange()
            \If{cooldown $>$ cooldownPeriod}
                \State rand $\gets$ RandomUniform(0,1)
                \If{rand $<$ veh.getLaneChangeProbability}
                    \State newLane $\gets$ getRandomLane()
                    \State roomInLane $\gets$ checkRoomInLane()
                    \If{roomInLane is True}
                        \State safetyCheck $\gets$ checkSafetyCondition()
                        \State incentiveCheck $\gets$ checkIncentiveCondition()
                        \If{safetyCheck \textbf{and} incentiveCheck}
                            \State successfulLaneChange $\gets$ True
                        \EndIf
                    \EndIf
                \EndIf
            \EndIf
            \State \Return successfulLaneChange
        \EndProcedure
    \end{algorithmic}
\end{algorithm}

\subsection{Model Parameters and Experimental Results}
Model parameters for Bando-FtL are estimated from existing work \cite{blondin_feedback_2019}. As described previously, the two populations of vehicles differ in two factors: (1) vehicle length and (2) $\alpha, \beta$ parameters within the velocity term of Bando-FtL. (1) The ``Car'' class of vehicles has a length of 4.5 meters, whereas the ``Truck'' class of vehicles has a length of 13.6 meters. (2) the $\alpha, \beta$ parameters for the ``Car'' class of vehicle are 0.5 and 20, respectively (as in \cite{blondin_feedback_2019}); we set $\alpha, \beta$ to values of 0.25 and 10, respectively for the ``Truck'' class of vehicles to approximate the slower acceleration. Additionally, the $v_{max}$ parameter within the optimal velocity term (Equation \ref{eq:opt-velocity}) is 5 and 4.5 for cars and trucks, respectively -- this decision was made to approximate the slower traveling speed of trucks compared to cars \cite{razzaq_model_2017, fu_empirical_2020, mehar_speed_2013}.

To investigate the relationship between velocity variation and the proportion of trucks within traffic, we simulate three scenarios of different penetration rates that roughly fall within the extreme values for observable truck penetration rates \cite{razzaq_model_2017}. For each scenario, there are 100 total vehicles: (1) 10\% trucks with each lane being 750 meters, (2) 20\% trucks with each lane scaled to 895 meters, and (3) 30\% trucks with each lane scaled to 1040 meters. We scale the road length proportionately to the total vehicle length to preserve the ratio of available road to total vehicle length. For each simulation, we track two summary metrics: (1) the maximum velocity variation of each vehicle, computed as the maximum velocity minus the minimum velocity, and (2) the total variation of the velocity, computed as the absolute value of the difference at each time point of the velocity. We simulate each scenario 100 times, with randomly initialized vehicle positions and noisy velocities, and aggregate the two summary metrics across the trial runs within scenarios. 

\begin{figure}[!htb]
    \centering
    \includegraphics[width=0.8\linewidth, trim={5cm 0cm 7cm 0cm}, clip]{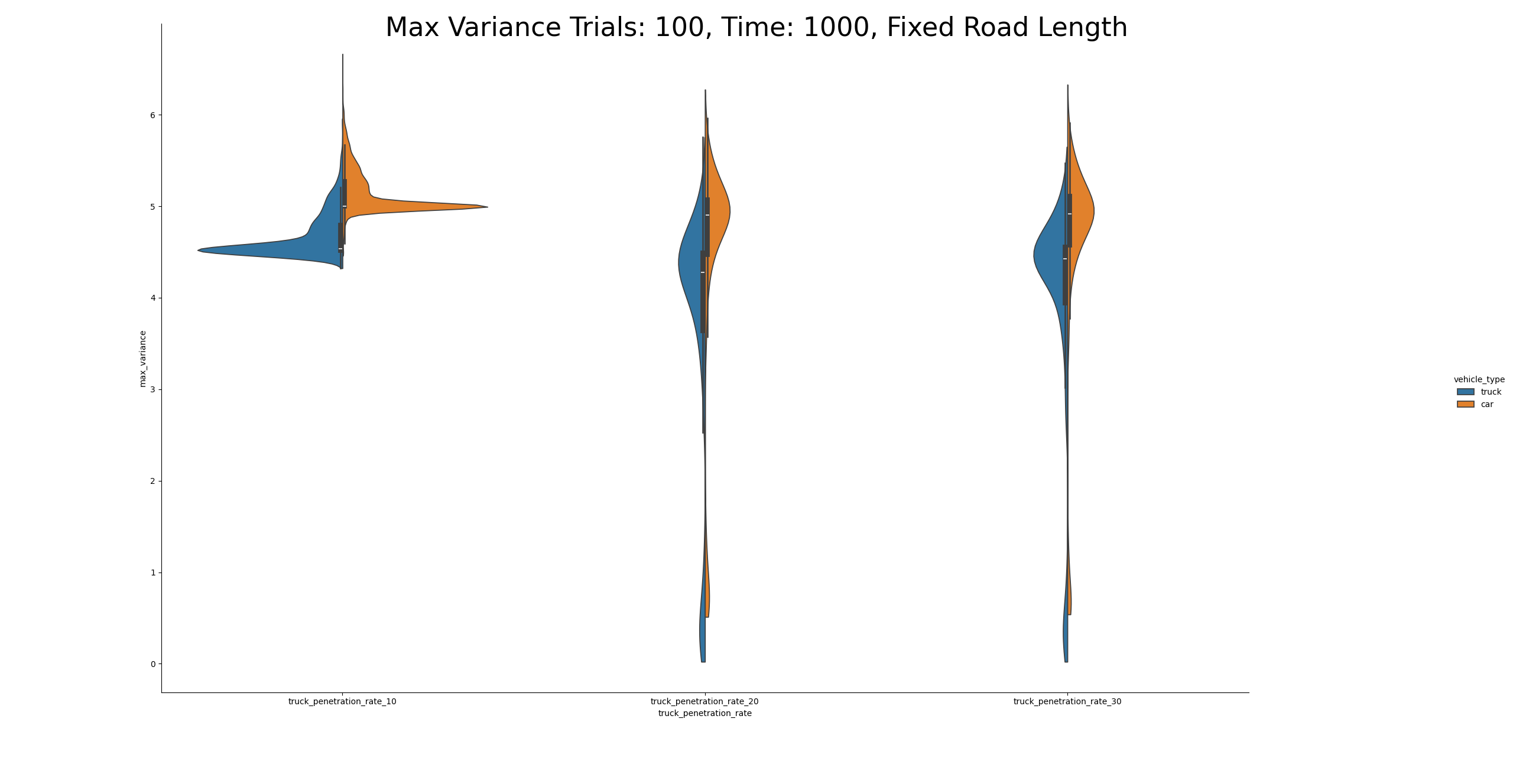}
    \caption{Branched violin plots showing the difference between maximum velocity and minimum velocity for each vehicle across 100 trials of each scenario. Violin plots correspond to 10\%, 20\%, and 30\% truck penetration rates from left to right. Each violin plot shows the max velocity difference for cars (blue, left side) and trucks (orange, right side).}
    \label{fig:max-variance}
\end{figure}

\begin{figure}[!htb]
    \centering
    \includegraphics[width=0.8\linewidth, trim={7cm 0cm 7cm 0cm}, clip]{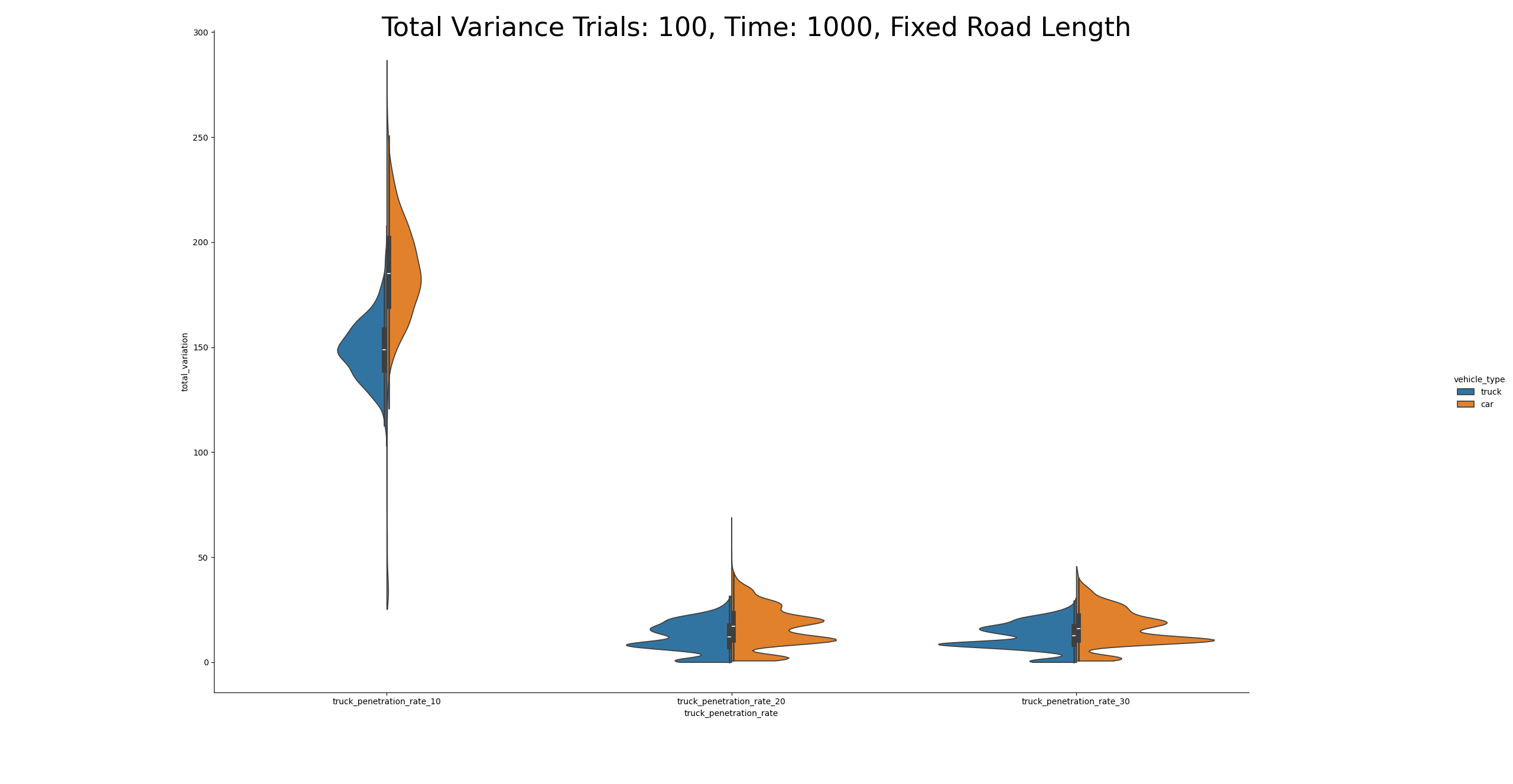}
    \caption{Branched violin plots showing the total variation of velocity for each vehicle across 100 trials of each scenario. Violin plots correspond to 10\%, 20\%, and 30\% truck penetration rates from left to right. Each violin plot shows the total variation of velocity for cars (blue, left side) and trucks (orange, right side).}
    \label{fig:total-variation}
\end{figure}

Figure \ref{fig:max-variance} reports the maximum difference in velocity across the different scenarios. We note that the maximum difference in velocity has the tightest distribution in the 10\% truck scenario; additionally, as the truck penetration rate increases from 20\% to 30\%, we observe a tighter distribution in the maximum difference in velocity for the ``car'' class of vehicles. This result suggests that a larger number of cars saw less variance in their speed when measuring by maximum velocity difference. Figure \ref{fig:total-variation} reports the total variation in velocity across the different scenarios, aggregating across the 100 individual trials of each scenario. We note that the largest total variation (for both cars and trucks) occurs during the 10\% truck scenario; we observe that the cars show a lower total variation in velocity (lower mean and tighter distribution) when compared to the trucks, which can be intuitively explained by the fact that cars can lane-change more readily (they require less gap in an adjacent lane) and that the incentive condition promote lane-changing when acceleration in a different lane would be higher. Additionally, the total variation for velocity is similar for both cars and trucks in the 20\% and 30\% truck penetration rate scenarios.


 \section{Conclusion and future work}

In this paper, we analyzed the existence of an optimal control for a minimization problem associated with the dynamics of a multi-population, multi-lane traffic model, both at the microscopic scale and in the macroscopic (mean-field) limit. This work represents a preliminary step toward a deeper understanding of the mechanisms governing multi-population traffic flow on multi-lane roads and how such systems can be improved.

Several questions remain open. For instance, identifying appropriate parameter ranges in the convolution kernels, as well as in the Incentive \eqref{lane-changing_acc} and Safety \eqref{lane-changing_safe} conditions, is still an open challenge. Addressing this issue requires comparison with empirical data. Furthermore, the probability function \eqref{eqn: probability_function} could be enhanced by incorporating context-dependent parameters, such as the vehicle’s position along the road, proximity to intersections, weather conditions, and time of day. These additions would allow for the implementation of more effective optimal control strategies and the development of predictive numerical tools.

Another promising direction for future research is the extension of the Mean Field Pontryagin Maximum Principle in the spirit of \cite{bongini2017mean} to the hybrid setting considered in this work. Such an extension would provide a deeper theoretical foundation for characterizing optimal controls in systems that combine continuous dynamics with discrete events.

Finally, we emphasize that controlling traffic in the presence of large vehicles is critically important not only for efficiency but also for safety. Indeed, studies and data indicate that heavy vehicles often increase both the likelihood and severity of road accidents \cite{KHORASHADI2005910}.

\bibliographystyle{plain}
\bibliography{reference.bib}
\end{document}